\documentclass[11pt]{article}

\usepackage{graphicx}

\usepackage{labelfig}

\usepackage{amsmath,amsfonts,amssymb,amsthm,upref,bm}%

\usepackage{mathtools}

\usepackage{mathrsfs}

\usepackage{booktabs}

\usepackage[paperwidth=210mm,paperheight=297mm,top=20mm,bottom=20mm,left=20mm,right=20mm]{geometry}

\DeclareMathAlphabet{\varmathbb}{U}{pxsyb}{m}{n}
\DeclareMathAlphabet{\mathpzc}{OT1}{pzc}{m}{it}

\newcommand{\MF}[1]{\mathop{#1}\nolimits}

\newcommand{\ii}{\kern0.05em\mathrm{i}\kern0.05em} \newcommand{\E}[1]{\textrm{e}^{#1}}

\renewcommand{\vec}[1]{\bm{#1}}

\renewcommand{\Re}{\MF{\mathrm{Re}}} \renewcommand{\Im}{\MF{\mathrm{Im}}}

\newcommand{\sign}{\mathop{\textrm{sign}}}

\newcommand{\Zz}{\varmathbb{Z}}%
\newcommand{\Cb}{\varmathbb{C}}%

\newcommand{\Pcoef}{\mathcal{P}}
\newcommand{\Qcoef}{\mathcal{Q}}
\newcommand{\Rcoef}{\mathcal{R}}
\newcommand{\Scoef}{\mathcal{S}}

\newcommand{\SectInf}{S}%

\newcommand{\BCzero}{\mathscr{B}_{0\infty}}%
\newcommand{\BCone}{\mathscr{B}_{1\infty}}%
\newcommand{\BCa}{\mathscr{B}_{a\infty}}

\def\Heun{H}%
\def\HeunL{\MF{\textit{Hl\/}}}%
\def\HeunS{\MF{\textit{Hs\/}}}%
\def\psHeunL{\mathpzc{Hl}_{\!0}}%
\def\psHeunLfin{\mathpzc{Hl}}%
\def\psHeunSfin{\mathpzc{Hs}}%
\def\psHeun{\mathpzc{H}}%
\def\psHeunLnearone{\mathpzc{Hl}^{(1)}}%
\def\psHeunSnearone{\mathpzc{Hs}^{(1)}}%
\def\psHeunLneara{\mathpzc{Hl}^{(a)}}%
\def\psHeunSneara{\mathpzc{Hs}^{(a)}}%
\def\psHeunLnearinf{\mathpzc{Hl}^{(\infty)}}%
\def\psHeunSnearinf{\mathpzc{Hs}^{(\infty)}}%
\def\midpoint#1{z^\star_{#1}}%
\def\bz{\tilde{b}}%
\def\nsp{n_*}%
\def\Lkr{\mathscr{L}}%
\def\reservecoeff{\varkappa}
\def\car{\mathpzc{c}}
\def\caz{\skew4\tilde{\car}}%

\allowdisplaybreaks[3]

\begin{document}

\title{On evaluation of the Heun functions}

\author{Oleg\ V.\ Motygin}
 
\date{Institute of Problems in Mechanical Engineering, Russian Academy\\ of
Sciences,
 V.O., Bol'shoy pr., 61, 199178 St.\,Petersburg, Russia\\ 
email: \texttt{o.v.motygin@gmail.com}}

\maketitle

\begin{abstract}
In the paper we deal with the Heun functions\,---\,solutions of the Heun
equation, which is the most general Fuchsian equation of second order with four
regular singular points. Despite the increasing interest to the equation and
numerous applications of the functions in a wide variety of physical problems,
it is only Maple amidst known software packages which is able to evaluate the
Heun functions numerically. But the Maple routine is known to be imperfect: even
at regular points it may return infinities or end up with no result. Improving
the situation is difficult because the code is not publicly available. The
purpose of the work is to suggest and develop alternative algorithms for
numerical evaluation of the Heun functions. A procedure based on power series
expansions and analytic continuation is suggested  which allows us to avoid
numerical integration of the differential equation and to achieve reasonable
efficiency and accuracy. Results of numerical tests are given.
\end{abstract}

\section{Introduction}
\label{intro0}

In the present paper we deal with the Heun functions which are solutions of the
equation introduced by Karl Heun in 1889 \cite{Heun} as a generalization of the
hypergeometric equation. Heun's equation is the most general Fuchsian equation
of second order with four regular singular points; we refer to
\cite{Ronveaux1995,SlavyanovLay2000,192-Maier2007,SleemanKuznetzov2010} for a
comprehensive mathematical treatment of the topic. At the same time, it is only
in recent years when the equation has become popular in the physics literature.
Now the Heun equation appears in many fields of modern physics and it is used to
describe a wide variety of physical phenomena\,---\,a comprehensive literature
on physical applications can be found in \cite{Hortacsu2011}. As a good source
of information on the current development of the field one should also mention
``The Heun project'': http://theheunproject.org/.

Despite the increasing interest to the equation, at the moment the only, to
author's knowledge, software package able to evaluate the Heun functions
numerically is Maple. However, the code is far from being perfect, it is not
difficult to encounter problems, when for ordinary parameters, at regular points
the routine returns infinities or spends tens seconds with no result. The quality
of the code is a serious obstacle for (potentially very numerous) applications
of the functions. Improving the situation is virtually impossible because the
code is not publicly available. 

The purpose of the present work is to develop alternative algorithms for
numerical evaluation of the Heun functions. We suggest a procedure based on
power series expansions and analytic continuation which allows us to avoid
numerical integration of the differential equation and to achieve reasonable
efficiency and accuracy. Program code is presented in \cite{mycode}. Results of
numerical tests and comparison with a case when Heun's function reduces to a
simple algebraic function are given.

The algorithm is applicable for computation of the multi-valued Heun functions.
In the present paper we also define their single-valued counterparts by fixation
of branch cuts. (Notably, in most studies of the Heun functions the subject is
not discussed; the author has also been unable to find information on branch
cuts in Maple's documentation on the Heun function.) For the single-valued
functions an improvement of the algorithm for points close to the singular ones
is described.

It should be noted that the developed algorithms are not intended to be
universal\,---\,more or less ordinary parameters are assumed. Surely, numerical
problems are expected and special treatment is needed for the cases of merging
singular points (see \cite{LaySlavyanov1999}) or large accessory parameter.

\section{Statement and basic notations} 
\label{sect:statement}

We start by writing Heun's equation in the standard form (see \cite{Erdelyi1955})
\begin{equation}
  \MF{\Heun''}(z)+\left(\frac{\gamma}{z}+\frac{\delta}{z-1}+
  \frac{\varepsilon}{z-a}\right)\MF{\Heun'}(z)+
  \frac{\alpha\beta z-q}{z(z-1)(z-a)}\MF{\Heun}(z)=0
\label{eq:Heun}
\end{equation}
with the Riemann $P$-symbol:
\[
P\left\{\begin{matrix}0 & 1 & a & \infty\\
0 & 0 & 0 & \alpha\\
1-\gamma & 1-\delta & 1-\varepsilon & \beta\\
\end{matrix};\ \ \ z\right\}.
\]

The parameter $q\in\Cb$ is usually referred to as an accessory or
auxiliary parameter and $\alpha$, $\beta$, $\gamma$, $\delta$, $\varepsilon$
(also belonging to $\Cb$) are exponent-related parameters connected via the
Fuchsian relation
\begin{equation*}
\alpha+\beta+1=\gamma+\delta+\varepsilon.
\end{equation*}
The equation has four regular singular points located at $z=0$, $1$, $a$,
$\infty$. It will be assumed below that $a \in\Cb$, $a\neq\{0,1,\infty\}$. In
the notation of the solutions sometimes we will omit some of the parameters so
that $\Heun(z)=\Heun(a,\ldots;z)=\Heun(a,q,\alpha,\beta,\gamma,\delta;z)$.

The Frobenius method can be used to derive local power-series solutions to
\eqref{eq:Heun}. There are $8$ local solutions of equation \eqref{eq:Heun} (two
per a singular point). In \S\;\ref{sect:expat0} we will present the two local
power-series solutions in a neighbourhood of the point $z = 0$. One of them is
analytic in a vicinity of zero and if $\gamma$ is not a nonpositive integer,
this solution, normalized to unity, is called the local Heun function (see
\cite{Ronveaux1995}). It is usually denoted by
$\HeunL(a,q,\alpha,\beta,\gamma,\delta;z)$. For the second Frobenius local
solution we will use the notation $\HeunS(a,q,\alpha,\beta,\gamma,\delta;z)$.

When $\gamma$ is a nonpositive integer, one solution of \eqref{eq:Heun} is
analytic but equal to zero at $z=0$, whereas the second solution can be
normalized to unity but generally is not analytic. It can be an arguable point,
but we will denote by $\HeunL(a,q,\alpha,\beta,\gamma,\delta;z)$ the normalized
solution and by $\HeunS(a,q,\alpha,\beta,\gamma,\delta;z)$ the analytic one.

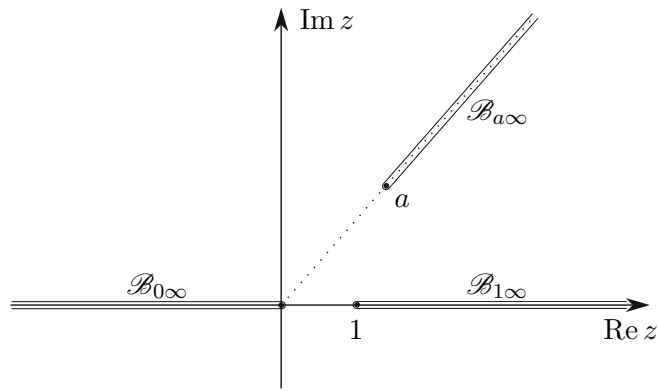
\begin{figure}[b]
\centering
 \SetLabels
 \L (0.53*0.145) $1$\\
 \L (0.6*0.48) $a$\\
 \L (0.455*0.93) $\Im z$\\
 \L (0.92*0.145) $\Re z$\\
 \L (0.19*0.26) $\BCzero$\\
 \L (0.71*0.26) $\BCone$\\
 \L (0.71*0.7) $\BCa$\\
 \endSetLabels
 \leavevmode\AffixLabels{\includegraphics[width=86mm]{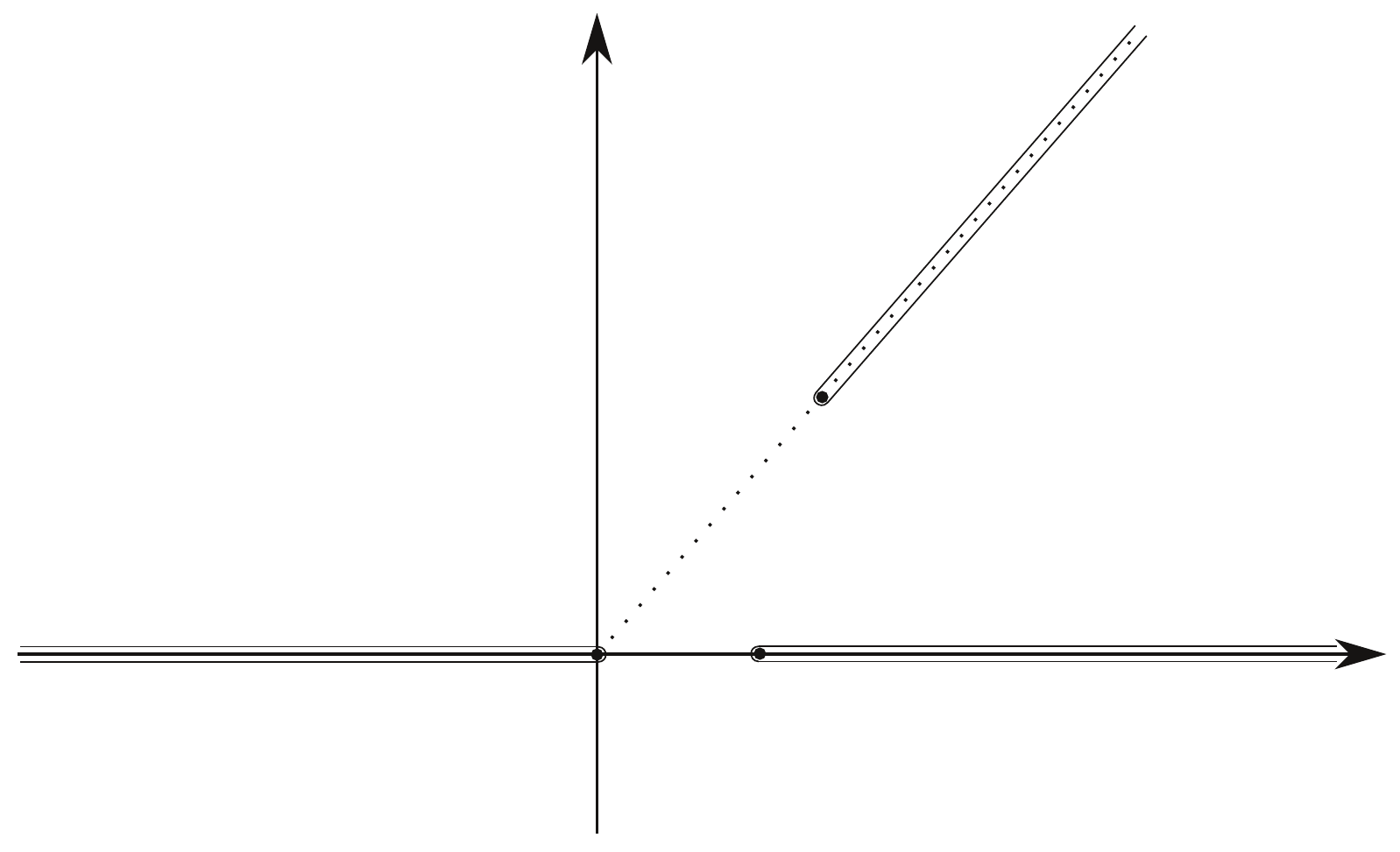}}
 \vspace{-1mm}
 \caption{Branch cuts.}
 \vspace{0mm}
 \label{bcwithoutlog}
\end{figure}

It is known that generally $\HeunL(a,\ldots;z)$ is a multi-valued function with
branch points at $z=1$, $a$ and $\infty$ and, so, to define a single-valued
function we should choose branch cuts. In the present work we fix the branch
cuts $\BCone=(1,+\infty)$ and $\BCa=(a,\E{\ii\arg(a)}\infty)$ connecting the
points $1$ and $a$ to $\infty$, respectively (see fig.~\ref{bcwithoutlog}). The
second function $\HeunS(a,\ldots;z)$ has the same branch cuts $\BCone$ and
$\BCa$. Besides, the function $\HeunS(z)$\,---\,generally, and the function
$\HeunL(z)$\,---\,for $\gamma\in\{0\}\cup\Zz^-$ ($\Zz^-$ means the set of
negative integers), have a branch point at $z=0$. So, we will also define the
branch cut $\BCzero=(-\infty,0)$.

It should be noted that the choice of branch cuts in Maple is unclear, but,
definitely, differs from that in the present paper. With the present fixation
the function definition domain is star-like with respect to zero which is
natural as we will define functions $\HeunL(z)$, $\HeunS(z)$ by analytic
continuation from a vicinity of $z=0$ (see \S\,\ref{sect:basic_algorithm}).
Besides, there is good agreement between branch cuts of the single-valued
functions $\HeunL$, $\HeunS$ and their representations near singular points (see
\S\,\ref{sect:improvement}).

\section{Power series expansion at the point $\bm{z=0}$}
\label{sect:expat0}

Power series expansion of the local Heun function $\HeunL$, such that
$\HeunL(a,q,\alpha,\beta,\gamma,\delta;0)=1$, is well-known since \cite{Heun}
for $\gamma\not\in\{0\}\cup\Zz^-$. We have
\begin{equation}
\HeunL(a,q,\alpha,\beta,\gamma,\delta;z)=\sum_{n=0}^{\infty}b_n z^n,
\label{eq:HeunL0series}
\end{equation}
where the coefficients $b_n$ are subject to the following three-term
recurrence:
\begin{equation}
  P_n b_n=Q_n b_{n-1}+R_n b_{n-2}.
\label{eq:HeunL0Andef}
\end{equation}
Here
\begin{equation}
\begin{gathered}
P_n = a n (\gamma-1+n),\\ Q_n = q +
(n-1)\bigl[(a+1)(\gamma+n-2)+\varepsilon+a\delta\bigr],\\ R_n =
-(n-2+\alpha)(n-2+\beta)
\end{gathered}
\label{eq:PQRdef}
\end{equation}
and the initial conditions are
\[
b_{-1}=0,\qquad b_0=1.
\]

The Heun function
$\HeunL(z)$ is analytic in the circle $|z|<R_0$, where $R_0$ is the distance
from zero to the nearest singular point, $R_0=\min\{1,|a|\}$. Then the famous
Cauchy's theorem on the expansion of an analytic function into a power series
(see e.g.\ Theorem 16.7 in \cite[Part~I]{Markushevich1977}) guarantees that the
series \eqref{eq:HeunL0series} converges to $\HeunL$ inside the circle
$|z|<R_0$. At this, of course, there might be issues with stability of the
recurrence process. For the stability it can be useful to write
\eqref{eq:HeunL0series} in the form $\HeunL(z)=\sum_{n=0}^{\infty}\bz_n(z)$,
where $\bz_n(z)=b_n z^n$ and the recurrence takes the form: $P_n \bz_n=z Q_n
\bz_{n-1}+z^2 R_n \bz_{n-2}$.

In the case $\gamma\in\Zz$ the local Frobenius solution corresponding to the
smaller exponent ($0$ or $1-\gamma$) may contain a logarithmic factor (see e.g.\
\cite{Kamke1944}). So, for $\gamma\in\{0\}\cup\Zz^-$ we are looking for the
solution of \eqref{eq:Heun}, satisfying $\HeunL(0)=1$, in the following form:
\begin{equation}
\HeunL(a,q,\alpha,\beta,\gamma,\delta;z)=\sum_{n=0,\,n\neq\nsp}^{\infty}\!\!c_n
z^n +\log(z)\sum_{n=\nsp}^{\infty}s_n z^n,
\label{eq:HeunLserieslog}
\end{equation}
where $\nsp=1-\gamma$. Note that a solution with the sought property
$\HeunL(0)=1$ could be found for any $c_{\nsp}$. We fix $c_{\nsp}=0$ for
definiteness.

Substituting \eqref{eq:HeunLserieslog} into \eqref{eq:Heun} we find
\begin{equation}
\log(z)\MF{\Lkr}\Biggl(\,\sum_{n=\nsp}^\infty s_n z^n\Biggr)+
\MF{\Lkr}\Biggl(\,\sum_{n=0,\,n\neq\nsp}^\infty\! c_n z^n\Biggr)+
\MF{\Hat{\Lkr}}\Biggl(\,\sum_{n=\nsp}^\infty s_n z^n\Biggr)=0,
\label{eq:HeunEqLserieslog}
\end{equation}
where $\Lkr$ is the operator of the Heun equation such that \eqref{eq:Heun} is
written as $\Lkr\Heun=0$. Besides,
\[
\left(\Hat{\Lkr}\psi\right)(z)=\frac2z\MF{\psi'}(z)+\frac1z\left(
\frac{\gamma-1}{z}+\frac{\delta}{z-1}+\frac{\varepsilon}{z-a}\right)\MF{\psi}(z).
\]
Obviously, \eqref{eq:HeunEqLserieslog} separates into two equations
\begin{gather}
\MF{\Lkr}\Biggl(\,\sum_{n=\nsp}^\infty s_n
z^n\Biggr)=0,\label{eq:HeunEqLserieslog1}\\
\MF{\Lkr}\Biggl(\,\sum_{n=0,\,n\neq\nsp}^\infty\! c_n z^n\Biggr)+
\MF{\Hat{\Lkr}}\Biggl(\,\sum_{n=\nsp}^\infty s_n z^n\Biggr)=0,
\label{eq:HeunEqLserieslog2}
\end{gather}

Now we collect in \eqref{eq:HeunEqLserieslog} terms having identical asymptotic
nature as $z\to0$. First we find that coefficients $c_n$ for $n=1,\ldots,\nsp-1$
are submitted to the recurrence \eqref{eq:HeunL0Andef}: $P_n c_n=Q_n c_{n-1}+R_n
c_{n-2}$, where $P_n$, $Q_n$, $R_n$ are defined by \eqref{eq:PQRdef} and the
initial conditions are $c_{-1}=0$, $c_0=1$.

Further we find from \eqref{eq:HeunEqLserieslog1} and
\eqref{eq:HeunEqLserieslog2} that the coefficients $s_n$ for
$n=\nsp+1,\nsp+2,\ldots$ are submitted to the same recurrence relationships
\eqref{eq:HeunL0Andef}: $P_n s_n=Q_n s_{n-1}+R_n s_{n-2}$, where $s_{\nsp-1}=0$
and another initial conditions include coefficients $c_{\nsp-1}$, $c_{\nsp-2}$:
\begin{gather*}
a\nsp\, s_{\nsp}=c_{\nsp-1}\bigl[q-\gamma(\varepsilon+a\delta-a-1)\bigr]-
c_{\nsp-2}\bigl[(1+\gamma)(2-\delta-\varepsilon)+\alpha\beta\bigr].
\end{gather*}

At the next step we can define coefficients $c_n$ for $n=\nsp+1,\nsp+2,\ldots$
From \eqref{eq:HeunEqLserieslog2} we obtain the following relationship:
\begin{equation}
P_n c_n=Q_n c_{n-1}+R_n c_{n-2}+ S_n s_{n}+T_n s_{n-1}+U_n s_{n-2},
\label{eq:recurrlog}
\end{equation}
where
\begin{gather*}
S_n = a(1-\gamma-2n),\qquad
T_n = \varepsilon+a\delta+(a+1)(\gamma+2n-3),\qquad
U_n = 4-2n-\alpha-\beta.
\end{gather*}

In this way for $\gamma\in\{0\}\cup\Zz^-$, using \eqref{eq:HeunLserieslog}  we
obtain a zero-exponent Frobenius solution, equal to unity at $z=0$. The second
solution locally can be defined as follows (see also
\eqref{eq:HeunEqLserieslog1}):
\begin{equation}
  \HeunS(a,q,\alpha,\beta,\gamma,\delta;z)=\sum_{n=\nsp}^{\infty}\check{s}_n z^n,
\label{eq:HeunSsergammanonpositive}
\end{equation}
where $P_n \check{s}_n=Q_n \check{s}_{n-1}+R_n \check{s}_{n-2}$ for $n>\nsp$ and
$\check{s}_{\nsp}=1$, $\check{s}_{\nsp-1}=0$. 

As it was mentioned above, $c_{\nsp}$ in \eqref{eq:HeunLserieslog} could be arbitrary. In other
words, the linear combination
\[
  \HeunL(a,q,\alpha,\beta,\gamma,\delta;z)+C\HeunS(a,q,\alpha,\beta,\gamma,\delta;z)
\]
is equal to unity at $z=0$ for an arbitrary constant $C$. In this sense, for
$\gamma\in\{0\}\cup\Zz^-$ the choice of solution $\HeunL$ is non-unique.

Consider now the function $\HeunS(z)$ for arbitrary $\gamma$. We
should distinguish two situations: $\gamma=1$ and $\gamma\neq1$. For
$\gamma\neq1$ we can use the following representation (see Table~2 in
\cite{192-Maier2007}, index $[0_-][1_+][a_+][\infty_-]$):
\begin{equation}
 \HeunS(a,q,\alpha,\beta,\gamma,\delta;z) = z^{1-\gamma}\HeunL\bigl(a,q-(\gamma-1)(\varepsilon+a\delta),
 \beta-\gamma+1,\alpha-\gamma+1,2-\gamma,\delta;z\bigr).
\label{eq:HeunSgammaneq1}
\end{equation}
Notably, the later formula includes \eqref{eq:HeunSsergammanonpositive} as a
particular case, justifying our way to introduce $\HeunL$ and $\HeunS$ for
non-positive integer $\gamma$.

For $\gamma=1$, repeating the arguments used to derive representation of
$\HeunL$ in the case $\gamma\in\{0\}\cup\Zz^-$, we can find the following local
representation
\begin{equation}
\HeunS(a,q,\alpha,\beta,\gamma,\delta;z)=\sum_{n=1}^{\infty}d_n z^n
+\log(z)\sum_{n=0}^{\infty}t_n z^n.
\label{eq:HeunSgamma=1ser}
\end{equation}
Here $P_n t_n=Q_n t_{n-1}+R_n t_{n-2}$, where $t_{-1}=0$, $t_{0}=1$ and (cf.\
\eqref{eq:recurrlog})
\begin{equation*}
P_n d_n=Q_n d_{n-1}+R_n d_{n-2}+ S_n t_{n}+T_n t_{n-1}+U_n t_{n-2},\quad
d_{-1}=d_0=0.
\end{equation*}

\section{Power series expansion at an arbitrary point}
\label{sect:expatarbpt}

Further we will extend the local Heun functions outside the circle of
convergence of the series \eqref{eq:HeunL0series}, \eqref{eq:HeunLserieslog},
\eqref{eq:HeunSgamma=1ser}. To avoid numerical integration of the differential
equation \eqref{eq:Heun} we will use analytic continuation process and the power
series expansion which is derived in this section.

We seek the solution $\Heun_{(z_0,H_0,H'_0)}(a,q,\alpha,\beta,\gamma,\delta;z)$
to the equation \eqref{eq:Heun} satisfying the conditions
\begin{equation}
 \Heun(a,q,\alpha,\beta,\gamma,\delta;z_0)=H_0,\quad
 \frac{\partial}{\partial z}\Heun(a,q,\alpha,\beta,\gamma,\delta;z)\Bigr|_{z=z_0}=H'_0.
\label{eq:z0H0H'0}
\end{equation}
Here $z_0$ is an arbitrary point that is assumed not to belong to the set
$\{0,a,1,\infty\}$. We will derive power series expansion of the Heun function in
the following form:
\begin{equation}
\Heun_{(z_0,H_0,H'_0)}(a,q,\alpha,\beta,\gamma,\delta;z)=\sum_{n=0}^{\infty}\car_n(z-z_0)^n.
\label{eq:HeunZ0series}
\end{equation}   

Substituting \eqref{eq:HeunZ0series} into \eqref{eq:Heun} and collecting terms
at the same power of $z-z_0$, we obtain the following 4-term recurrent
relationship defining $\car_n$:
\begin{equation}
\Pcoef_n \car_n=\Qcoef_n \car_{n-1}+\Rcoef_n \car_{n-2} + \Scoef_n
\car_{n-3},\label{eq:HeunZ0Cndef}
\end{equation}
where
\begin{gather*}
\Pcoef_n=-n(n-1)z_0(z_0-1)(z_0-a),\\ \Qcoef_n=(n-1)\Bigl\{
\bigl[\gamma+\delta+\varepsilon+3(n-2)\bigr]z_0^2+
\bigl[(a+1)(4-2n-\gamma)-\varepsilon-a\delta\bigr]z_0 + a(\gamma+n-2)\Bigr\}, \\
\Rcoef_n=\bigl\{(n-2)\bigl[2(\gamma+\delta+\varepsilon)+3(n-3)\bigr]+
\alpha\beta\big\}z_0-q
{}-(n-2)\bigl[(a+1)(\gamma+n-3)+\varepsilon+a\delta\bigr], \\
\Scoef_n=(n-3)\bigl(\gamma+\delta+\varepsilon+n-4\bigr)+\alpha\beta.
\end{gather*}
Obviously, the solution \eqref{eq:HeunZ0series} satisfies the conditions \eqref{eq:z0H0H'0} if
the recurrence process starts with the initial conditions 
\begin{equation*}
\car_{-1}=0,\qquad \car_{0}=H_0,\qquad \car_{1}=H'_0.
\end{equation*}

The series \eqref{eq:HeunZ0series} converges inside the circle $|z-z_0|<R$,
where $R$ is the distance to the nearest singular point,
$R=\min\{|z_0|,|z_0-1|,|z_0-a|\}$.

\section{Basic algorithm}
\label{sect:basic_algorithm}

First we consider evaluation of $\HeunL(z)$ for $\gamma\not\in\{0\}\cup\Zz^-$.
We introduce the projection operator $\mathscr{P}_{z}^N$ which, being applied to
an analytic function, truncates its power series expansion at the point $z$ to
the first $N$ terms. Using the expansion \eqref{eq:HeunL0series}, we evaluate
\begin{equation}
\bigl(\mathscr{P}_0^N\!\HeunL\bigr)(z)=\sum_{n=0}^{N}b_n z^n,\quad
\bigl(\mathscr{P}_0^N\!\HeunL\bigr)'(z)=\sum_{n=1}^{N}n\, b_n z^{n-1},
\label{eq:HeunLpwrserappr}
\end{equation}
as  approximation of $\HeunL(z)$ and $\HeunL'(z)$  in a vicinity of $z=0$.

In our algorithm the number $N$ in the representations
\eqref{eq:HeunLpwrserappr} is not fixed, it will be defined as we proceed with
recurrent computation of $\bz_n(z)=b_n z^n$ and summation until a termination
condition is satisfied. Namely, we stop the process when
$\bigl(\mathscr{P}_0^{N}\!\HeunL\bigr)'(z)$ and
$\bigl(\mathscr{P}_0^{N-1}\!\HeunL\bigr)'(z)$ are not distinguishable in the
used computer arithmetics and $|\bz_n(z)|<\epsilon$, where $\epsilon$ is the
machine epsilon.

To estimate the quality of the approximation, in view of \eqref{eq:Heun} we
compute the value
\[
  \skew5\Hat{\HeunL}(z)=\frac{1}{q-\alpha\beta z}\biggl\{z(z-1)(z-a)\bigl(\mathscr{P}_0^N\!\HeunL\bigr)''(z)+
  \Bigl[\gamma(z-1)(z-a)+\delta z(z-a)+\varepsilon z(z-1)\Bigr]\bigl(\mathscr{P}_0^N\!\HeunL\bigr)'(z)\biggr\},
\]
where $\bigl(\mathscr{P}_0^N\!\HeunL\bigr)''(z)=\sum_{n=2}^{N}n(n-1)\,b_n z^{n-2}$.
Then we suppose proximity of
\begin{equation}
\MF{r_0}(z)=\bigl|\skew5\Hat{\HeunL}(z)-\bigl(\mathscr{P}_0^N\!\HeunL\bigr)(z)\bigr|
\label{eq:errdef}
\end{equation}
to the true error of the approximation
$\bigl|\HeunL(z)-\bigl(\mathscr{P}_0^N\!\HeunL\bigr)(z)\bigr|$. Note that
numerical computation of $\skew5\Hat{\HeunL}(z)$ near the point
$z=z_*=q/(\alpha\beta)$ is unreliable due to essential loss of significance. It
can be suggested for a vicinity of $z_*$, say, for $\{z:|q-\alpha\beta
z|<0.01\}$,  to use an estimate based on properties of the series, e.g.
\begin{equation}
\MF{\Hat{r}\vphantom{r}_0}(z)=
\sqrt{N}|\bz_N(z)| + \epsilon N \bigl|\bigl(\mathscr{P}_0^N\!\HeunL\bigr)(z)\bigr|.
\label{eq:err2def}
\end{equation}

We can write the described algorithm as a function $\psHeunL(z)$ which returns 4-tuple
\[
 \psHeunL: z\mapsto [f,f',r,N],
\]
where $N+1$ is the number of terms in power series, defined by the termination
condition, $r$ is the value computed with \eqref{eq:errdef} or
\eqref{eq:err2def}, $f=\bigl(\mathscr{P}_0^N\!\HeunL\bigr)(z)$ and
$f'=\bigl(\mathscr{P}_0^N\!\HeunL\bigr)'(z)$.

The scheme of computation of $\HeunL(z)$ in the case $\gamma\in\{0\}\cup\Zz^-$
is analogous, though a bit more involved. We use the expansion
\eqref{eq:HeunLserieslog} and, instead of \eqref{eq:HeunLpwrserappr}, define the
function $\psHeunL(z)$ starting from the expression
\[
\sum_{n=0,\,n\neq\nsp}^{N}\!\!c_n z^n +\log(z)\sum_{n=\nsp}^{N}s_n z^n.
\]

Assume that $|z|<\reservecoeff R_0$, where $R_0=\min\{1,|a|\}$ is the radius of
convergence of the series in \eqref{eq:HeunL0series} or
\eqref{eq:HeunLserieslog} and $\reservecoeff\in(0,1)$ is some coefficient chosen
so that $N$ defined by the termination condition can be expected to be moderate,
say $\reservecoeff=0.5$. Then we can use the numerical algorithm $\psHeunL$ for
evaluation of the function $\HeunL(z)$, its derivative  and for estimation of
the approximation error.

Consider further the case $|z|\geq\reservecoeff R_0$. First we define an
auxiliary algorithm. Let $z_0$ be an arbitrary point not belonging to
$\{0,1,a,\infty\}$. Using \eqref{eq:HeunZ0series} we define
\[
\bigl(\mathscr{P}_{\!z_0}^N\Heun_{(z_0,H_0,H'_0)}\bigr)(z)=\sum_{n=0}^{N}\car_n(z-z_0)^n,\quad
\bigl(\mathscr{P}_{\!z_0}^N\Heun_{(z_0,H_0,H'_0)}\bigr)'(z)=\sum_{n=1}^{N}n\,\car_n(z-z_0)^{n-1}.
\]
as approximations of $\Heun_{(z_0,H_0,H'_0)}(z)$ and
$\Heun'_{(z_0,H_0,H'_0)}(z)$ for $z$ close to $z_0$. Here coefficients $\car_n$
are defined by the recurrence relationships \eqref{eq:HeunZ0Cndef}. So we
proceed with recurrent computation of $\caz_n(z)=\car_n (z-z_0)^n$ and summation
until the termination condition (analogous to that described above) is
satisfied.

We compute
\begin{multline*}
  \skew4\Hat{\Heun}_{(z_0,H_0,H'_0)}(z)=\frac{1}{q-\alpha\beta z}\biggl\{z(z-1)(z-a)\bigl(\mathscr{P}_{\!z_0}^N\Heun_{(z_0,H_0,H'_0)}\bigr)''(z)\\
  {}+
  \Bigl[\gamma(z-1)(z-a)+\delta z(z-a)+\varepsilon z(z-1)\Bigr]\bigl(\mathscr{P}_{\!z_0}^N
  \Heun_{(z_0,H_0,H'_0)}\bigr)'(z)\biggr\},
\end{multline*}
and the value
\begin{equation}
\MF{r_{(z_0,H_0,H'_0)}}(z)=\bigl|\skew4\Hat{\Heun}_{(z_0,H_0,H'_0)}(z)-
\bigl(\mathscr{P}_{\!z_0}^N\Heun_{(z_0,H_0,H'_0)}\bigr)(z)\bigr|,
\label{eq:errz0def}
\end{equation}
supposing its proximity to the true error of the approximation
$\bigl|\Heun_{(z_0,H_0,H'_0)}(z)-
\bigl(\mathscr{P}_{\!z_0}^N\Heun_{(z_0,H_0,H'_0)}\bigr)(z)\bigr|$. In view of
essential loss of significance in computation of
$\skew4\Hat{\Heun}_{(z_0,H_0,H'_0)}(z)$, near $z=z_*$ we define
\begin{equation}
\MF{\Hat{r}\vphantom{r}_{(z_0,H_0,H'_0)}}(z)=
\sqrt{N}|\caz_N(z)| + \epsilon N \bigl|\bigl(\mathscr{P}_{\!z_0}^N\Heun_{(z_0,H_0,H'_0)}\bigr)(z)\bigr|.
\label{eq:err2z0def}
\end{equation}

We can write the described algorithm as a function $\psHeun_{(z_0,H_0,H'_0)}(z)$ which returns 4-tuple
\[
 \psHeun_{(z_0,H_0,H'_0)}: z\mapsto [f,f',r,N],
\]
where $N+1$ is the number of terms in power series defined by the termination
condition, $r$ is the value computed with \eqref{eq:errz0def} or
\eqref{eq:err2z0def},
$f=\bigl(\mathscr{P}_{z_0}^N\Heun_{(z_0,H_0,H'_0)}\bigr)(z)$ and
$f'=\bigl(\mathscr{P}_{z_0}^N\Heun_{(z_0,H_0,H'_0)}\bigr)'(z)$.

\begin{figure}[t]
\centering
 \SetLabels
 \L (0.57*0.305) $1$\\
 \L (0.41*0.64) $a$\\
 \L (0.3*0.95) $\Im z$\\
 \L (0.92*0.29) $\Re z$\\
 \L (0.24*0.305) $0$\\
 \L (0.7775*0.615) $z$\\
 \L (0.37*0.385) $z_1$\\
 \L (0.425*0.415) $z_2$\\
 \L (0.475*0.445) $z_3$\\
 \L (0.525*0.474) $z_4$\\
 \L (0.578*0.503) $z_5$\\
 \L (0.64*0.54) $z_6$\\
 \L (0.72*0.587) $z_7$\\
 \endSetLabels
 \leavevmode\AffixLabels{\includegraphics[width=120mm]{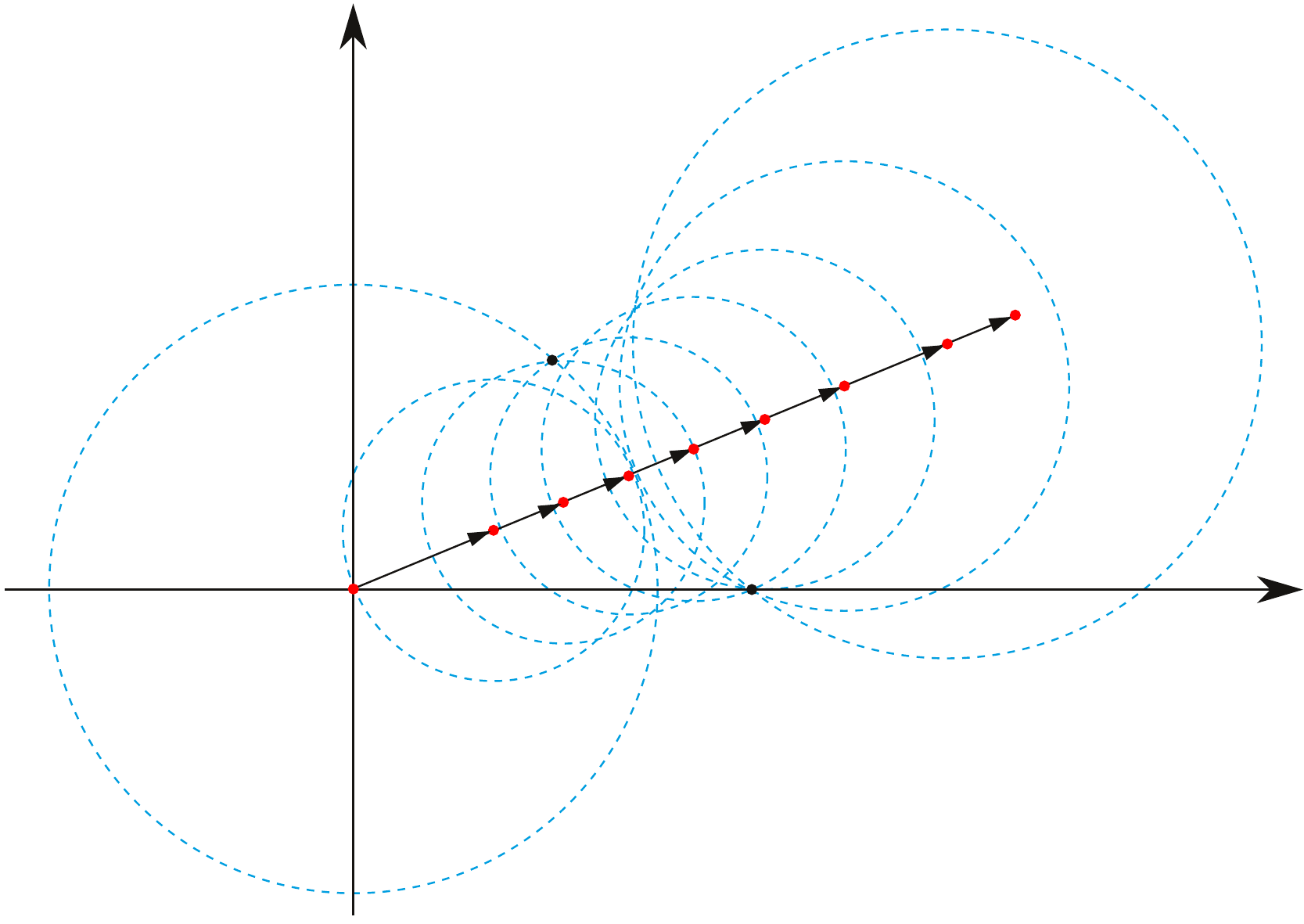}}
 \vspace{-1mm}
 \caption{Analytic continuation using power series.}
 \vspace{0mm}
 \label{analytcont}
\end{figure}
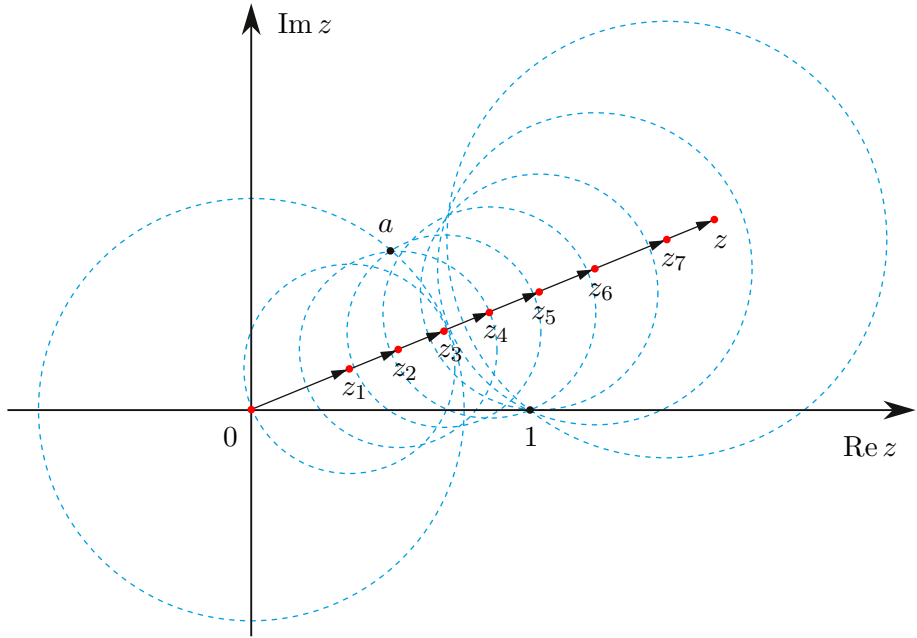


Now we are in a position to describe the algorithm based on analytic
continuation along a path from zero to $z$. Consider first the simplest version
of the algorithm when the path is the line segment $(0,z)$. At the first step,
we compute
\[
 [\HeunL_1,\HeunL'_1,r_1,N_1]=\psHeunL(z_1),
\]
where $z_1=\E{\ii\arg(z)}\reservecoeff R_0$ (see Fig.~\ref{analytcont}). 
Further, for $p=1$, $2$, and so on, we define 
\[
R_p=\min\{|z_p|,|z_p-1|,|z_p-a|\},
\]
\[
z_{p+1}=\begin{cases}z_p+\E{\ii\arg(z)}\reservecoeff R_p & \mbox{if}\ \ |z-z_p|>\reservecoeff R_p,\\
z & \mbox{if}\ \ |z-z_p|\leq\reservecoeff R_p,
\end{cases}
\]
and compute
\[
 [\HeunL_{p+1},\HeunL'_{p+1},r_{p+1},N_{p+1}]=\psHeun_{(z_p,\HeunL_p,\HeunL'_p)}(z_{p+1}).
\]
The iterations stops when $z_{p+1}=z$. Finally, we have $\HeunL_{p+1}$,
$\HeunL'_{p+1}$ as approximations of $\HeunL(z)$ and $\HeunL'(z)$, respectively.
We also compute the values $r_\Sigma=r_1+\ldots+r_{p+1}$ and
$N_\Sigma=N_1+\ldots+N_{p+1}+p+1$. Here $N_\Sigma$ is the total number of power
series terms which can be used as a measure of computer load and $r_\Sigma$ may
be an indicator of the approximation quality.


The described above algorithm of continuation along a line segment is readily
generalized for the case when $0$ and $z$ are connected by a polyline
$\Upsilon$. This gives us a way to compute the multi-valued Heun function. The
resulting procedure can be considered as a function
$\mathring{\psHeunLfin}(\Upsilon)$ which returns 4-tuple
\[
 \mathring{\psHeunLfin}: \Upsilon\mapsto [f,f',r_\Sigma,N_\Sigma],
\]
where $f$ and $f'$ are the resulting approximations of the Heun function at $z$
and its derivative. It should be noted that the function $\HeunS(z)$ for
$\gamma\neq1$ is expressed through $\HeunL(z)$ via \eqref{eq:HeunSgammaneq1} and
to define $\mathring{\psHeunSfin}(\Upsilon)$ through
$\mathring{\psHeunLfin}(\Upsilon)$ one should only compute the multi-valued
function $z^{-\gamma}$ along $\Upsilon$. In the case $\gamma=1$ the procedure of
analytic continuation described above for $\HeunL$ can be applied with simple
modification\,---\,it should start from the expansion
\eqref{eq:HeunSgamma=1ser}.

It is easy to observe that the size of the step in the described analytic
continuation is small for points of the polyline $\Upsilon$ close to one of the
singular points. This also means an increase of the number of used circular
elements in the continuation procedure which, in its turn, may lead to loss of
accuracy. The influence of the singular points can be reduced by a choice of the
path of continuation. 

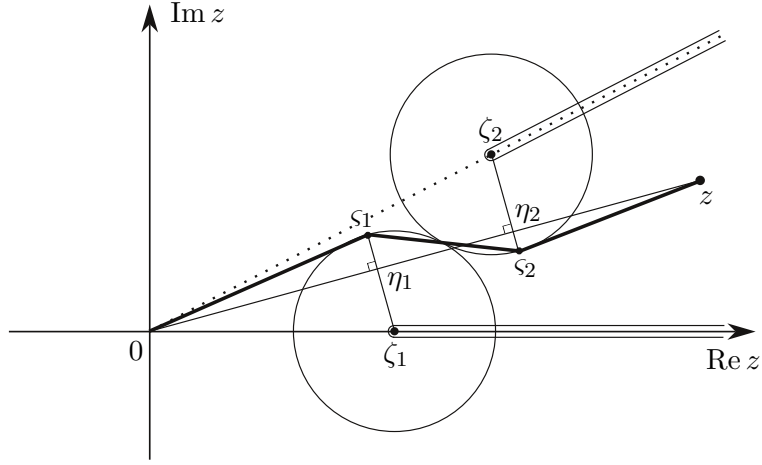
\begin{figure}[t]
\centering
 \SetLabels
 \L (0.215*0.935) $\Im z$\\
 \L (0.92*0.2) $\Re z$\\
 \L (0.495*0.215) $\zeta_1$\\
 \L (0.62*0.693) $\zeta_2$\\
 \L (0.16*0.22) $0$\\
 \L (0.91*0.55) $z$\\
 \L (0.5*0.378) $\eta_1$\\
 \L (0.668*0.523) $\eta_2$\\
 \L (0.452*0.51) $\varsigma_1$\\
 \L (0.667*0.405) $\varsigma_2$\\
 \endSetLabels
 \leavevmode\AffixLabels{\includegraphics[width=100mm]{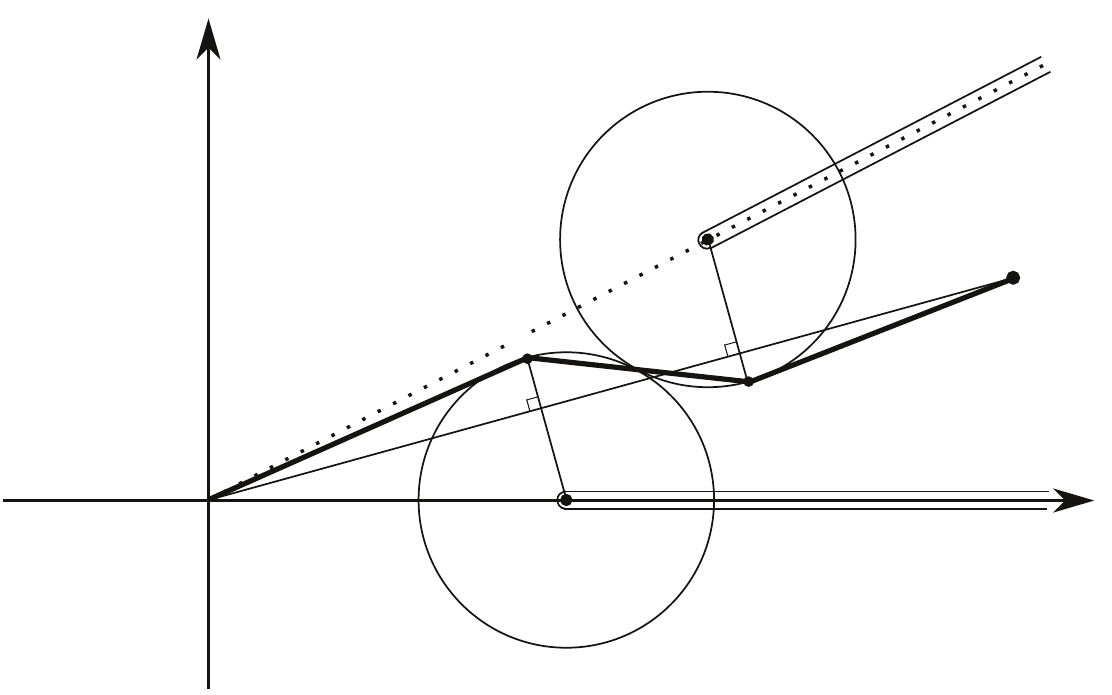}}
 \vspace{-1mm}
 \caption{Path from zero to $z$ consisting of three line segments}
 \vspace{0mm}
 \label{analytcont2}
\end{figure}

Let us give some details of application of the above algorithm for computation
of the single-valued Heun functions. We will only use simple paths consisting of
two or three line segments. We denote $\zeta_1=a$ and $\zeta_2=1$ if $|a|<1$ and
$\zeta_1=1$ and $\zeta_2=a$ otherwise. Let $R_j =
\min\{|\zeta_j|,|\zeta_1-\zeta_2|\}$ ($j=1,2$) and $\eta_j$ be the closest to
$\zeta_j$ point of the line segment $(0,z)$ (see fig.~\ref{analytcont2}). Note
that $\zeta_j\neq\eta_j$ due to the assumption $z\not\in\BCone\cup\BCa$. Let
$d_j=|\zeta_j-\eta_j|$ be the distance from $\zeta_j$ to $(0,z)$. If $\eta_j$ is
an internal point of $(0,z)$ and $d_j<R_j/2$, then we introduce new point
\[
\varsigma_j=\zeta_j+\exp\bigl(\ii\pi/2+\ii\arg(z)\bigr)\min\{R_j/2,|z-\zeta_j|\}
\sign(\Im(z/\zeta_j)).
\]
In this way we obtain the sequence of points $\Pi$, which is either $[0,z]$ or
$[0,\varsigma_1,z]$ or $[0,\varsigma_2,z]$ or $[0,\varsigma_1,\varsigma_2,z]$.
Then we define the path $\Upsilon$ that consequently connects the points of
$\Pi$ by line segments (see an example of three-segment path in
fig.~\ref{analytcont2}) and define
$\psHeunLfin(z)=\mathring{\psHeunLfin}(\Upsilon)$ and
$\psHeunSfin(z)=\mathring{\psHeunSfin}(\Upsilon)$.

\section{Computation of single-valued Heun functions near singular points}
\label{sect:improvement}

In this section we consider computation of single-valued Heun functions. We have
already noted that the number of circular elements in the analytic continuation
procedure becomes large as $z$ approaches one of the singular points
$\{1,a,\infty\}$. So, we suggest an improvement of the algorithm; for this we
will compute local solutions in a vicinity of the singular point, using results
of \cite{192-Maier2007} and writing these solutions in terms of $\HeunL(z)$ and
$\HeunS(z)$.

Let $z$ be close to $1$. We use the following representation of the Heun
function:
\begin{multline}
\HeunL(a,q,\alpha,\beta,\gamma,\delta;z) = C_1
\HeunL(1-a,\alpha\beta-q,\alpha,\beta,\delta,\gamma;1-z)+ C_2
\HeunS(1-a,\alpha\beta-q,\alpha,\beta,\delta,\gamma;1-z),
\label{eq:near1}
\end{multline}
where $C_1$, $C_2$ are some constants, the Heun functions in the right-hand side
(the first and the second local solutions in a vicinity of $z=1$) correspond to
the index $[1_+0_+][a_+][\infty_+]$ in Table~2 \cite{192-Maier2007}. There is a
complication when the point $1$ belongs to $\BCa$ ($a\in(0,1)$). In this case
the coefficients in \eqref{eq:near1} are generally different for $z$ belonging
to the upper and the lower half-spaces: $C^{(\pm)}_1$, $C^{(\pm)}_2$ for
$\{z:\pm\Im(z)>0\}$.

Consider branch cuts of $\HeunL(1-a,\ldots;1-z)$ and $\HeunS(1-a,\ldots;1-z)$.
The points of the branch cut, corresponding to $\BCa$ for $\HeunL(a,\ldots;z)$
and $\HeunS(a,\ldots;z)$, are defined by the equation $1-z=s(1-a)$, where
$s\in(1,+\infty)$. The set of points $z=1+s(a-1)$ constitutes the ray emanating
from $a$ to $\exp\{\ii\arg(a-1)\}\infty$ and, thus, it is located outside the
circle of radius $|a-1|$ with center at $z=1$. Analogously, consider the branch
cut of $\HeunL(1-a,\ldots;1-z)$ and $\HeunS(1-a,\ldots;1-z)$ corresponding to
the branch cut $\BCone$. The points of the branch cut are defined by the
equation $1-z=s$, where $s\in(1,+\infty)$, thus, lying outside the circle of
radius $1$ with center at $z=1$. The last possible branch cut of
$\HeunL(1-a,\ldots;1-z)$ and $\HeunS(1-a,\ldots;1-z)$ is defined as $1-z=q$,
$q\in(-\infty,0)$ and coincides with the branch cut $\BCone$ of the function
$\HeunL(z)$ in the left hand side of \eqref{eq:near1}.

In order to use the representation \eqref{eq:near1} we should define the
coefficients $C_1$, $C_2$ (or, in the special case, $C^{(\pm)}_1$,
$C^{(\pm)}_2$). Since for the general Heun equation an explicit solution to the
two-point connection problem is not known, we will find the matching
coefficients in the following way. First, we define a matching point
$z=\midpoint{1}$, preferably being equally close to zero and $1$ and distant
from $z=a$. In our algorithm we fix
\[
\midpoint{1} = \frac12 + \frac{\ii s}{\sqrt{2}},\ \ \mbox{where}\ \
s=\begin{cases}\phantom{-}\sign(\Im(z)) & \textrm{if}\ a\in(0,1),\\
 -\sign(\Im(a)) &
\textrm{otherwise}.\end{cases}
\]

Further, we apply the algorithms $\psHeunLfin$ and $\psHeunSfin$ described in
\S\;\ref{sect:basic_algorithm} and find
\begin{gather*}
 [f_0,f'_0,r_0,N_0]=\psHeunLfin(a,q,\alpha,\beta,\gamma,\delta,\midpoint{1}),\\
 [f_1,f'_1,r_1,N_1]=\psHeunLfin(1-a,\alpha\beta-q,\alpha,\beta,\delta,\gamma,1-\midpoint{1}),\\ 
 [f_2,f'_2,r_2,N_2]=\psHeunSfin(1-a,\alpha\beta-q,\alpha,\beta,\delta,\gamma,1-\midpoint{1}). 
\end{gather*}
Then we solve the linear system
\[
\left(\begin{matrix} f_1 & f_2\\
 -f'_1 & -f'_2
\end{matrix}\right)
\left(\begin{matrix}
C_1\\C_2
\end{matrix}\right)=
\left(\begin{matrix}
f_0\\f'_0
\end{matrix}\right)
\]
and obtain the matching coefficients $C_1$, $C_2$. It may also be reasonable to
keep the found values $C_1=C_1(a,q,\alpha,\beta,\gamma,\delta)$,
$C_2=C_2(a,q,\alpha,\beta,\gamma,\delta)$ in computer memory.

On finding $C_1$, $C_2$ (by computation or in the computer memory) we define
\begin{equation*}
 \psHeunLnearone: z\mapsto [f,f',r,N],
\end{equation*}
where 
\begin{gather*}
 f=C_1 f_1 + C_2 f_2,\quad 
 f'=-C_1 f'_1 - C_2 f'_2,\quad
 r = |C_1|r_1 + |C_2|r_2,\quad
 N = N_1 + N_2,\\
 [f_1,f'_1,r_1,N_1]=\psHeunLfin(1-a,\alpha\beta-q,\alpha,\beta,\delta,\gamma,1-z),\\
 [f_2,f'_2,r_2,N_2]=\psHeunSfin(1-a,\alpha\beta-q,\alpha,\beta,\delta,\gamma,1-z).
\end{gather*}

The described scheme can be repeated literally to define $\psHeunSnearone$ based
on the representation
\begin{equation}
\HeunS(a,q,\alpha,\beta,\gamma,\delta;z) = C'_1
\HeunL(1-a,\alpha\beta-q,\alpha,\beta,\delta,\gamma;1-z)+ C'_2
\HeunS(1-a,\alpha\beta-q,\alpha,\beta,\delta,\gamma;1-z),
\label{eq:near1s}
\end{equation}
where $C'_1$ and $C'_2$ are some coefficients to be found (or two sets
${C'}^{(-)}_1$, ${C'}^{(-)}_2$ and ${C'}^{(+)}_1$, ${C'}^{(+)}_2$\,---\,in the
special case $a\in(0,1)$).

We note that finding $C_1$, $C_2$ or $C'_1$, $C'_2$ demands computation of
$\HeunL(a,\ldots;\midpoint{1})$ or $\HeunS(a,\ldots;\midpoint{1})$, and both
terms $\HeunL(1-a,\ldots;1-\midpoint{1})$ and
$\HeunS(1-a,\ldots;1-\midpoint{1})$ in the right hand side of \eqref{eq:near1}
or \eqref{eq:near1s}. So, if the matching constants are not known, the
algorithms $\psHeunLnearone$ and $\psHeunSnearone$ are preferable over
$\psHeunLfin$ and $\psHeunSfin$ in a sufficiently small vicinity of $z=1$. In
the code used in \S\,\ref{sect:num} the algorithms are applied for $|z-1|<0.05\,R_1$
($R_1=\min\{1,|a-1|\}$) when the matching constants are not known and for
$|z-1|<0.25\,R_1$ when the coefficients are already precomputed.

In a very similar way we consider the case when $z$ is located near $a$. We will
use the following representation of the Heun function:
\begin{multline}
\HeunL(a,q,\alpha,\beta,\gamma,\delta;z)\\{} = D_1
 \HeunL\Bigl(\frac{a-1}{a},\alpha\beta-\frac{q}{a},
 \alpha,\beta,\varepsilon,\gamma;\frac{a-z}{a}\Bigr)
 +D_2\HeunS\Bigl(\frac{a-1}{a},\alpha\beta-\frac{q}{a},\alpha,\beta,\varepsilon,\gamma;\frac{a-z}{a}\Bigr),
\label{eq:neara}
\end{multline}
where $D_1$, $D_2$ are some constants and the two functions in the right-hand
side correspond to index $[a_+0_+1_+][\infty_+]$ in Table~2
\cite{192-Maier2007}. There is a complication when the point $a$ belongs to
$\BCone$ or $\BCzero$. In this case we generally have different coefficients in
\eqref{eq:near1} in the upper and the lower half-spaces: $D^{(\pm)}_1$,
$D^{(\pm)}_2$ for $\{z:\pm\Im(z)>0\}$.

Consider branch cuts of $\HeunL((a-1)/a,\ldots;\zeta)$ and
$\HeunS((a-1)/a,\ldots;\zeta)$, where $\zeta=(a-z)/a$. The points of the branch
cut that emanates from the singularity $\zeta=(a-1)/a$ are defined in the
$z$-plane by the equation $a-z=s(a-1)$, where $s\in(1,+\infty)$. This set of
points constitutes the ray emanating from $1$ to $\exp\{\ii\arg(1-a)\}\infty$
and located outside the circle of the radius $|1-a|$ with center at $z=a$. The
branch cut of $\HeunL((a-1)/a,\ldots;\zeta)$ and $\HeunS((a-1)/a,\ldots;\zeta)$
emanating from $\zeta=1$ is a ray going from zero to $\exp\{\ii\arg(-a)\}\infty$
($z=a(1-s)$); it is located outside the circle of radius $|a|$ with center at
$z=a$. The third possible branch cut of $\HeunL((a-1)/a,\ldots;\zeta)$ and
$\HeunS((a-1)/a,\ldots;\zeta)$ is defined as $a-z=aq$, $q\in(-\infty,0)$ and
coincides with the branch cut $\BCa$ of the function $\HeunL(a,\ldots;z)$.

In order to find the coefficients $D_1$, $D_2$ (or, in the special case,
$D^{(\pm)}_1$, $D^{(\pm)}_2$), first, we define a matching point
$z=\midpoint{a}$, preferably being equally close to zero and $a$ and distant
from $z=1$. We fix
\[
\midpoint{a} = \frac{a}{2} + \frac{\ii s}{\sqrt{2}},\ \ \mbox{where}\ \
s=\begin{cases}\sign(\Im(z)) & \textrm{if}\ a\in\BCzero\cup\BCone,\\
 a|a|^{-1}\sign(\Im(a)) &
\textrm{otherwise}.\end{cases}
\]

Further, we find
\begin{gather*}
 [f_0,f'_0,r_0,N_0]=\psHeunLfin(a,q,\alpha,\beta,\gamma,\delta,\midpoint{a}),\\
 [f_1,f'_1,r_1,N_1]=\psHeunLfin\Bigl(\frac{a-1}{a},\alpha\beta-\frac{q}{a},\alpha,\beta,
   \varepsilon,\gamma,\frac{a-\midpoint{a}}{a}\Bigr),\\
 [f_2,f'_2,r_2,N_2]=\psHeunSfin\Bigl(\frac{a-1}{a},\alpha\beta-\frac{q}{a},\alpha,\beta,
   \varepsilon,\gamma,\frac{a-\midpoint{a}}{a}\Bigr),
\end{gather*}
and $D_1$, $D_2$ are defined as solution to the linear system
\[
\left(\begin{matrix} f_1 & f_2\\ -a^{-1}f'_1 & a^{-1}f'_2
\end{matrix}\right)
\left(\begin{matrix} D_1\\D_2
\end{matrix}\right)=
\left(\begin{matrix} f_0\\f'_0
\end{matrix}\right).
\]
It may also be reasonable to keep the found values
$D_1=D_1(a,q,\alpha,\beta,\gamma,\delta)$,
$D_2=D_2(a,q,\alpha,\beta,\gamma,\delta)$ in computer memory.

On finding $D_1$, $D_2$ (by computation or in the computer memory) we define
\begin{equation*}
 \psHeunLneara: z\mapsto [f,f',r,N],
\end{equation*}
where 
\begin{gather*}
 f=D_1 f_1 + D_2 f_2,\quad 
 f'=-D_1 a^{-1}f'_1 - D_2a^{-1}f'_2,\quad
 r = |D_1|r_1 + |D_2|r_2,\quad
 N = N_1 + N_2,\\
 [f_1,f'_1,r_1,N_1]=\psHeunLfin\Bigl(\frac{a-1}{a},\alpha\beta-\frac{q}{a},\alpha,\beta,\varepsilon,\gamma,\frac{a-z}{a}\Bigr),\\
 [f_2,f'_2,r_2,N_2]=\psHeunSfin\Bigl(\frac{a-1}{a},\alpha\beta-\frac{q}{a},\alpha,\beta,\varepsilon,\gamma,\frac{a-z}{a}\Bigr).
\end{gather*}
The described scheme can be repeated literally to define $\psHeunSneara$ based
on the representation
\begin{multline*}
\HeunS(a,q,\alpha,\beta,\gamma,\delta;z)\\{} = D'_1
 \HeunL\Bigl(\frac{a-1}{a},\alpha\beta-\frac{q}{a},
 \alpha,\beta,\varepsilon,\gamma;\frac{a-z}{a}\Bigr)
 +D'_2\HeunS\Bigl(\frac{a-1}{a},\alpha\beta-\frac{q}{a},\alpha,\beta,\varepsilon,\gamma;\frac{a-z}{a}\Bigr),
\end{multline*}
where $D'_1$ and $D'_2$  (or, in the special case, ${D'}^{(\pm)}_1$,
${D'}^{(\pm)}_2$) are some coefficients to be found.

We note again that the algorithms $\psHeunLneara$ and $\psHeunSneara$ are
preferable over $\psHeunLfin$ and $\psHeunSfin$ in a sufficiently small vicinity
of $z=a$. In the code used in \S\,\ref{sect:num} the algorithms are applied for
$|z-a|<0.05\,R_a$ ($R_a=\min\{|a|,|1-a|\}$) when the matching constants are not
known and for $|z-a|<0.25\,R_a$ when the coefficients are already precomputed.

Consider now vicinity of the point $z=\infty$. The branch cuts $\BCzero$,
$\BCone$ and $\BCa$ split the vicinity of infinity
$V_\infty=\{z:|z|>\max(1,|a|)\}$ into three sectors, further denoted as
$\SectInf^{(1)}$, $\SectInf^{(2)}$, $\SectInf^{(3)}$. Coefficients connecting
the function $\HeunL(z)$ with two local solutions at infinity are defined for
each of the sectors separately. Note that there is a special case when $a$ is
real, then the number of sectors decreases to two.

We use the following representation of the Heun function near $z=\infty$ for
$z\in\SectInf^{(k)}$ ($k=1, 2, 3$):
\begin{multline}
\HeunL(a,q,\alpha,\beta,\gamma,\delta;z) = E^{(k)}_1 z^{-\alpha}
\HeunL\Bigl(\frac{1}{a},\frac{q+\alpha(\delta-\beta)}{a}+
\alpha(\varepsilon-\beta),\alpha,\alpha-\gamma+1,\alpha-\beta+1,\delta;\frac{1}{z}\Bigr)
\\ {}+E^{(k)}_2 z^{-\alpha}
\HeunS\Bigl(\frac{1}{a},\frac{q+\alpha(\delta-\beta)}{a}+
\alpha(\varepsilon-\beta),\alpha,\alpha-\gamma+1,\alpha-\beta+1,\delta;\frac{1}{z}\Bigr),
\label{eq:nearinf}
\end{multline}
where $E^{(k)}_1$, $E^{(k)}_2$ are some constants. Further we will omit the
superscript for brevity. The two functions in the right-hand side correspond to
the index $[\infty_+0_+][1_+][a_+]$ in Table~2 \cite{192-Maier2007}.

Obviously, the branch cuts $\BCa$ and $\BCone$ transform for
$\HeunL(1/a,\ldots,1/z)$ and $\HeunS(1/a,\ldots,1/z)$ to line segments
connecting zero with $z=a$ and $z=1$ respectively. These branch cuts are located
outside $V_\infty$. The branch cut $\BCzero$ suits both sides of \eqref{eq:nearinf}.

To find the coefficients $E_1$, $E_2$, first, we define a matching point
$z=\midpoint{\infty}$. In our algorithm we fix
\[
\midpoint{\infty} = C_\infty R_\infty \E{\ii\omega_k}, \quad\mbox{for}\ \
z\in\SectInf^{(k)},\ \ k=1,2,3,
\]
where $\omega_k$ is the angle of the mean line of the sector $\SectInf^{(k)}$,
$R_\infty=\max\{1,|a|\}$ and $C_\infty>1$ ($C_\infty=2$ in the computations,
presented below). Then, we find
\begin{gather*}
 [f_0,f'_0,r_0,N_0]=\psHeunLfin(a,q,\alpha,\beta,\gamma,\delta,\midpoint{\infty}),\\
 [f_1,f'_1,r_1,N_1]=\psHeunLfin\biggl(\frac{1}{a},\frac{q+\alpha(\delta-\beta)}{a}+
\alpha(\varepsilon-\beta),\alpha,\alpha-\gamma+1,\alpha-\beta+1,\delta,\frac{1}{\midpoint{\infty}}\biggr),\\
 [f_2,f'_2,r_2,N_2]=\psHeunSfin\biggl(\frac{1}{a},\frac{q+\alpha(\delta-\beta)}{a}+
\alpha(\varepsilon-\beta),\alpha,\alpha-\gamma+1,\alpha-\beta+1,\delta,\frac{1}{\midpoint{\infty}}\biggr),
\end{gather*}
and the matching coefficients are defined as solution to the linear system
\begin{equation*}
\left(\begin{matrix} [\midpoint{\infty}]^{-\alpha}f_1 &
[\midpoint{\infty}]^{-\alpha}f_2\\[2mm] -[\midpoint{\infty}]^{-\alpha-1}
(f'_1/\midpoint{\infty}+\alpha f_1) ~&~ -[\midpoint{\infty}]^{-\alpha-1}
(f'_2/\midpoint{\infty}+\alpha f_2)
\end{matrix}\right)
\left(\begin{matrix} E_1\\[2mm]E_2
\end{matrix}\right)=
\left(\begin{matrix} f_0\\[2mm]f'_0
\end{matrix}\right).
\end{equation*}

Finally, we define a function $\psHeunLnearinf$ returning values of
$f(z)\approx\HeunL(z)$, $f'(z)\approx\HeunL'(z)$ for sufficiently large $|z|$.
On finding $E_1$, $E_2$ (by computation or in the computer memory) we compute
\begin{equation*}
 \psHeunLnearinf: z\mapsto [f,f',r,N],
\end{equation*}
where 
\begin{gather*}
 f= E_1 z^{-\alpha}f_1 + E_2 z^{-\alpha}f_2,\quad
 f'=
 -E_1 z^{-\alpha-1}(f'_1/z+\alpha f_1) - E_2 z^{-\alpha-1} (f'_2/z+\alpha f_2),\\
 r=\bigl|E_1 z^{-\alpha}\bigr|r_1 +\bigl|E_2 z^{-\alpha}\bigr|r_2,\quad
 N = N_1 + N_2,\\
 [f_1,f'_1,r_1,N_1]=\psHeunLfin\Bigl(\frac{1}{a},\frac{q+\alpha(\delta-\beta)}{a}+
\alpha(\varepsilon-\beta),\alpha,\alpha-\gamma+1,\alpha-\beta+1,\delta,\frac{1}{z}\Bigr),\\
 [f_2,f'_2,r_2,N_2]=\psHeunSfin\Bigl(\frac{1}{a},\frac{q+\alpha(\delta-\beta)}{a}+
\alpha(\varepsilon-\beta),\alpha,\alpha-\gamma+1,\alpha-\beta+1,\delta,\frac{1}{z}\Bigr).
\end{gather*}

The scheme can be repeated literally to define $\psHeunSnearinf$ based on the
representation for $z\in\SectInf^{(k)}$
\begin{multline*}
\HeunS(a,q,\alpha,\beta,\gamma,\delta;z) = {E'}^{(k)}_1 z^{-\alpha}
\HeunL\Bigl(\frac{1}{a},\frac{q+\alpha(\delta-\beta)}{a}+
\alpha(\varepsilon-\beta),\alpha,\alpha-\gamma+1,\alpha-\beta+1,\delta;\frac{1}{z}\Bigr)
\\ {}+{E'}^{(k)}_2 z^{-\alpha}
\HeunS\Bigl(\frac{1}{a},\frac{q+\alpha(\delta-\beta)}{a}+
\alpha(\varepsilon-\beta),\alpha,\alpha-\gamma+1,\alpha-\beta+1,\delta;\frac{1}{z}\Bigr),
\end{multline*}
where ${E'}^{(k)}_1$ and ${E'}^{(k)}_2$ are some coefficients to be found.

We note again that the algorithms $\psHeunLnearinf$ and $\psHeunSnearinf$ are
preferable over $\psHeunLfin$ and $\psHeunSfin$  for sufficiently large $|z|$.
In the code used in the following section the algorithms are applied for
$|z|>5\, C_\infty\,R_\infty$ when the matching coefficients are not known and
for $|z|>C_\infty\,R_\infty$ otherwise.

\section{Numerical results}
\label{sect:num}

In this section we present some results of numerical evaluation of the function
$\HeunL$. For this we will use the algorithm $\psHeunLfin$ with the
improvements described in the previous section ($\psHeunLnearone$,
$\psHeunLneara$ and $\psHeunLnearinf$).

It is known that in some cases Heun functions can be expressed in terms of
hypergeometric functions (see
\cite{Kuiken1979,ShaninCraster2002,Maier2005,Valent2007,Birkandan2014} and
references therein). For a test of the present algorithm we use the formula
given in \cite{Birkandan2014}:
\[
\HeunL\biggl(4,\frac94,\frac32,\frac32,\frac12,2,z\biggr)=
 {}_2F_1\biggl(\frac12,\frac12,\frac12,1-(z-1)^2\Bigl(1-\frac{z}{4}\Bigr)\biggr)
 =h(z):=\frac{2}{\sqrt{4-z}\,(1-z)},
\]
where ${}_2F_1$ is the ordinary hypergeometric function.

\begin{figure}[t!]
\centering
 \SetLabels
 \L (-0.06*0.87) $\Im z$\\
 \L (0.61*-0.03) $\Re z$\\
 \endSetLabels
 \leavevmode\AffixLabels{\includegraphics[width=120mm]{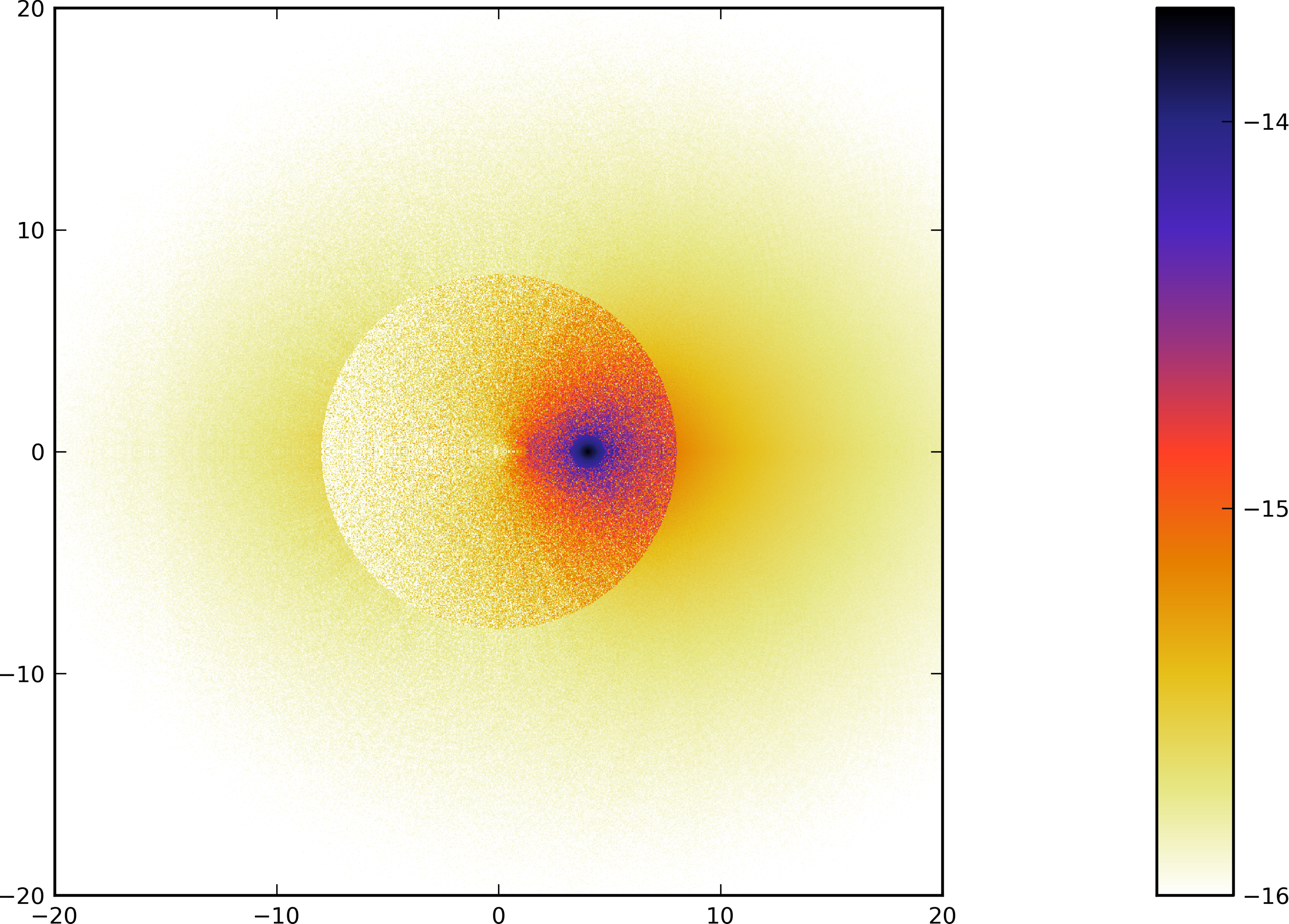}}
 \vspace{1mm}
 \caption{Values of $\max\bigl\{-16,\log_{10}\Lambda(z)\bigr\}$.}
 \vspace{2mm}
 \label{numcomp}
\end{figure}

\begin{figure}[t!]
\centering
 \SetLabels
 \L (-0.06*0.87) $\Im z$\\
 \L (0.61*-0.03) $\Re z$\\
 \endSetLabels
 \leavevmode\AffixLabels{\includegraphics[width=120mm]{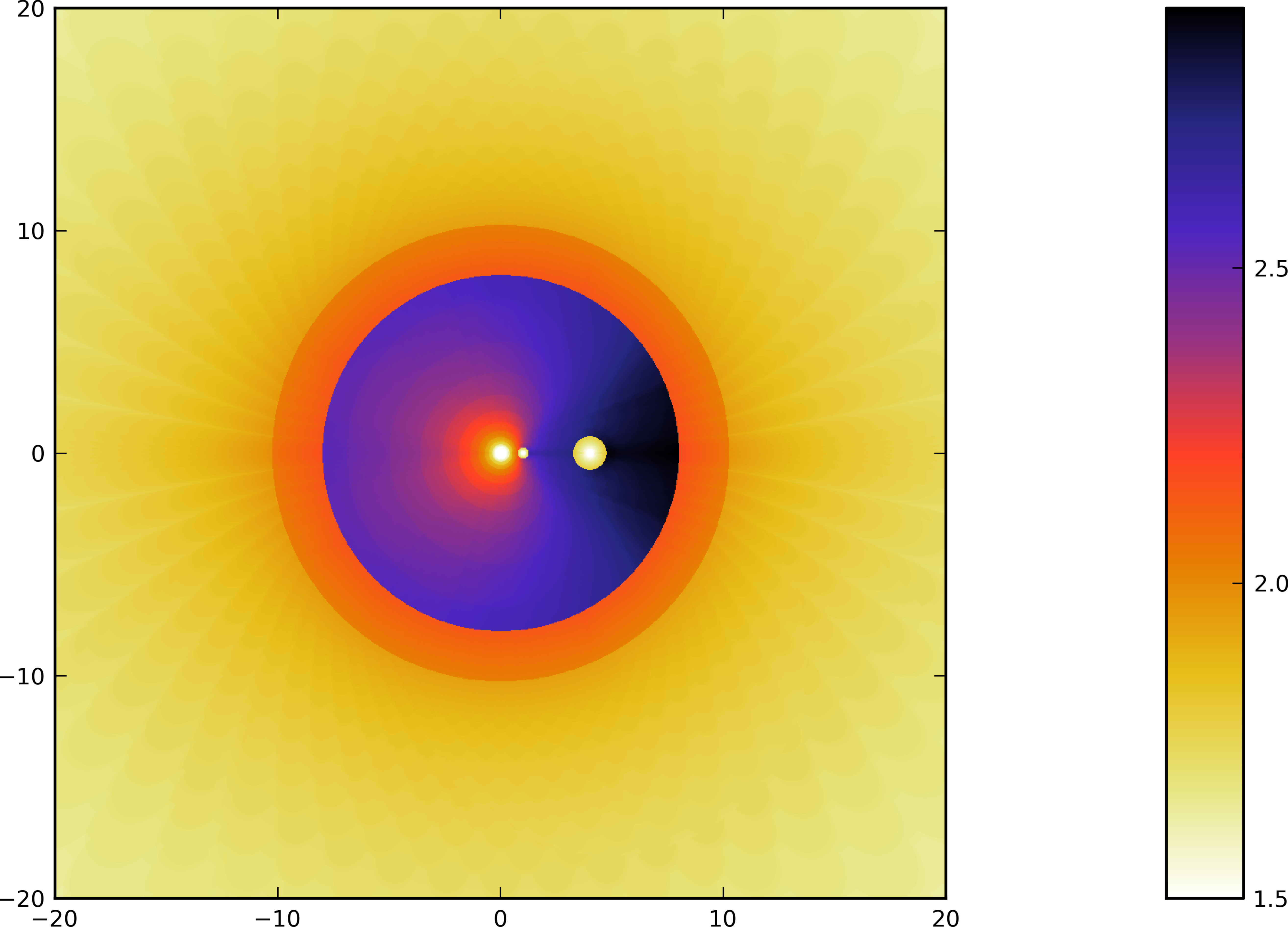}}
 \vspace{1mm}
 \caption{Values of $\max\bigl\{1.5,\log_{10}N(z)\bigr\}$.}
 \vspace{0mm}
 \label{numcomp3}
\end{figure}

We define
\[
\Delta(z):=\HeunL\biggl(4,\frac94,\frac32,\frac32,\frac12,2,z\biggr)-\frac{2}{\sqrt{4-z}\,(1-z)}
\]
and check numerically the identity $\Delta(z)=0$. Calculations are performed in
the numerical computing environment GNU Octave and double precision (64-bit)
arithmetics (the machine epsilon $\epsilon$ is about $2.22\cdot10^{-16}$).

Figure~\ref{numcomp} shows in a semilogarithmic scale results of computations
of the relative error
\[
 \Lambda(z)=\frac{|\Delta(z)|}{1+|h(z)|}+\frac{|\Delta'(z)|}{1+|h'(z)|}.
\]
on the grid $(\Re z,\Im z)
\in\vec{L}(1000,[-20,20])\times\vec{L}(1000,[-20,20])$, where $\vec{L}(n,\chi)$
is the set consisting of $n$ linearly spaced in the interval $\chi$ values
(including interval's end-points). It is easy to note that accuracy loss occurs
mainly near the singularity $z=4$. The maximum error for $\Lambda(z)$ found on
the grid is $1.9635\cdot10^{-14}$. So, as one can see, the accuracy of
computation of the function and its derivative is rather satisfactory for the
used double-float arithmetics.

In Fig.~\ref{numcomp3} we present (in a semilogarithmic scale) the value $N(z)$
which means the total number of terms in power series used to compute
$\HeunL\bigl(4,\frac94,\frac32,\frac32,\frac12,2,z\bigr)$. This value varies
from $N(0)=1$ and can characterize the time of computation. To be more specific,
in the Table we present typical times of evaluation in Octave of
$\HeunL\bigl(4,\frac94,\frac32,\frac32,\frac12,2,z\bigr)$ via the basic
algorithm $\psHeunLfin$ and via the algorithm $\psHeunLfin$ improved by
$\psHeunLnearone$, $\psHeunLneara$ and $\psHeunLnearinf$. The computer in use
has a 3.4\;GHz Intel Core i5 CPU and 8\;Gb of RAM.
\begin{table}[h]
\begin{center}
\begin{tabular}{c|c|c|c}
 $z$ & Time (sec.), basic algorithm & \multicolumn{2}{c}{Time (sec.), improved algorithm}    \\
  &  & 1st evaluation & subsequent evaluations  \\\midrule
 $20\ii$ & $0.04$ & $0.05$ & $0.007$\\
 $20+\ii\epsilon$ & $0.075$ & $0.05$ & $0.007$\\
 $-20$ & $0.035$ & $0.05$ & $0.003$\\
 $0.99$ & $0.03$ & $0.02$ & $0.003$\\
 $4+0.01\ii$ & $0.075$ & $0.05$ & $0.002$  
\end{tabular}
\end{center}
\end{table}

\vspace{-6mm}

\section*{Acknowledgements} The author is indebted to Professor A.M.\,Ishkhanyan
for drawing attention to the subject and useful discussions. The author thanks
Professor A.Ya.\,Kazakov for interest in the results of the present work
and useful remarks.

\end{document}